\documentclass[12pt]{amsart}

\usepackage{graphicx}
\usepackage{amsfonts}
\usepackage{epsf}
\usepackage{amssymb}
\usepackage{amsmath}
\usepackage{amscd}

\newtheorem{theorem}{Theorem}[section]
\newtheorem{proposition}[theorem]{Proposition}
\newtheorem{lemma}[theorem]{Lemma}

\newtheorem{corollary}[theorem]{Corollary}

\newtheorem{remark}[theorem]{Remark}

\newtheorem{obstruction}[theorem]{Sliceness obstruction}

\newcommand{\kerr}{\mbox{Ker} }
\newcommand{\cokerr}{\mbox{Coker} }

\newcommand{\s}{\mathfrak{s}}

\def\del{{\partial}}
\def\into{{\hookrightarrow}}
\def\im{{\textup{im}}}
\def\Tors{{\textup{Tors}}}
\def\coker{{\textup{coker}}}

\input xy
\xyoption{all}

\begin{document}
\title{The slice-ribbon conjecture for 3-stranded pretzel knots}
\begin{abstract}{We determine the (smooth) concordance order of the 3-stranded pretzel knots $P(p,q,r)$
with $p,q,r$ odd. We show that each one of finite order is, in fact, ribbon, thereby proving the
slice-ribbon conjecture for this family of knots. As corollaries we give new proofs of results first
obtained by Fintushel-Stern and Casson-Gordon. }
\end{abstract}
\author{Joshua Greene}
\address{Department of Mathematics\\ Princeton University\\ Princeton\\ NJ 08544}
\email{jegreene@math.princeton.edu}
\author{Stanislav Jabuka}
\address{Department of Mathematics and Statistics\\ University of Nevada\\ Reno\\ NV 89557}
\email{jabuka@unr.edu}
\maketitle
\section{Introduction}
\begin{figure}[htb!]
\centering
\includegraphics[width=10cm]{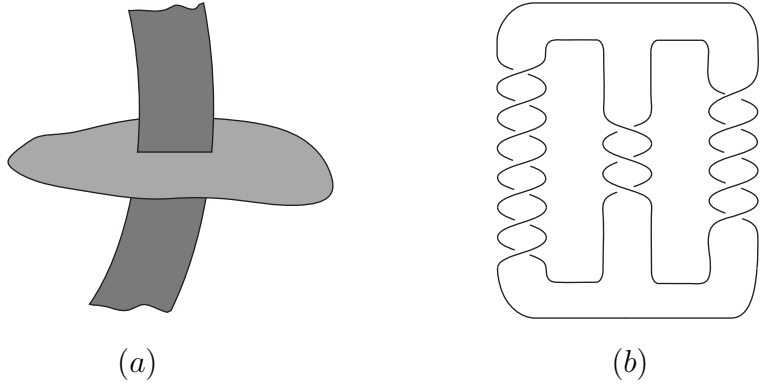}
\put(-243,-20){$(a)$}
\put(-56,-20){$(b)$}
\caption{(a) A ribbon singularity: An immersed disk $D^2 \looparrowright S^3$ intersects itself along an arc
of double points such that
one of the components of the preimage of that arc lies in the interior of the disk $D^2$.
(b) The pretzel knot $P(p,q,r)$ for the case of $p=7$, $q=-3$ and $r=5$.}  \label{pic2}
\end{figure}

Recall that a knot $K\subset S^3$ is called {\em slice} if it bounds a smoothly and properly
embedded disk $D^2 \hookrightarrow D^4$ in the 4-ball $D^4$, and
{\em ribbon} if it bounds an immersed disk $D^2 \looparrowright S^3$ with only ribbon singularities, as in
figure \ref{pic2}(a).
It is easy to see that every ribbon knot $K$ is
slice: simply push the ribbon singularities of the ribbon disk $D^2$ into the 4-ball to create an
embedded slice disk for $K$. The converse is not known to be true but is the content of an interesting
conjecture (Problem 1.33 on Kirby's list \cite{kirby}):
\vskip2mm
{\bf Slice-ribbon conjecture. } Every slice knot is ribbon.
\vskip2mm
In recent work \cite{lisca1, lisca2} Paolo Lisca proved the slice-ribbon conjecture for 2-bridge knots. A brief synopsis of his method of proof goes as follows:
given a 2-bridge knot $K$, find a
negative-definite plumbing description for $Y_K$, the 2-fold cover of $S^3$ branched over $K$ (which for
a 2-bridge knot is a lens space). If $K$ is slice then $Y_K$ also bounds a rational homology ball which,
when glued to the negative definite plumbing, yields a negative definite closed 4-manifold $X$. By
Donaldson's celebrated diagonalization theorem the intersection form of $X$ must be diagonalizable, an obstruction which
proves strong enough to single out all slice knots among 2-bridge knots. By finding explicit ribbons
for each of the several categories of slice knots, Lisca proves the slice-ribbon conjecture.  As pointed
out in \cite{lisca1}, for some 2-bridge knots $K$ it is necessary to consider the sliceness obstruction for
both $K$ and its mirror.

\begin{figure}[htb!]
\centering
\includegraphics[width=14cm]{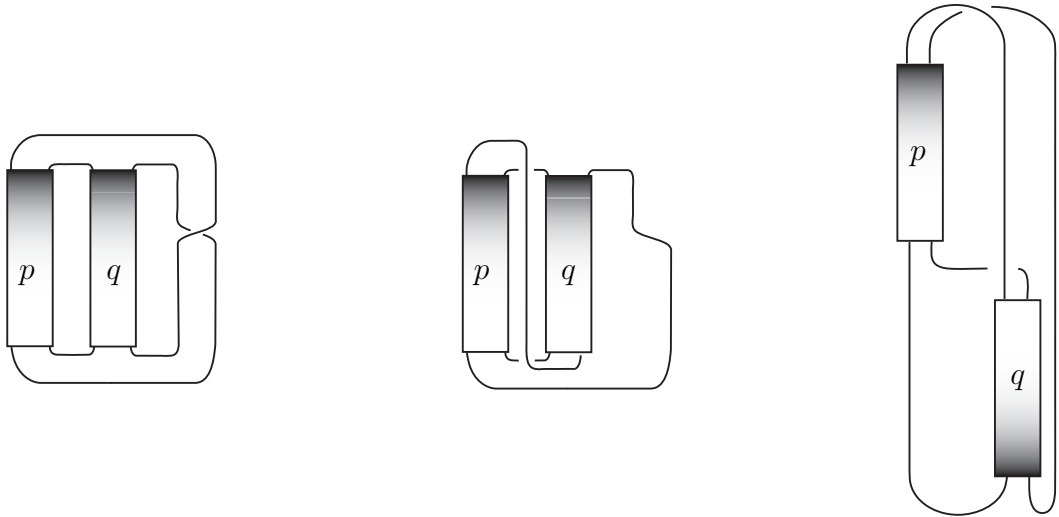}
\put(-393,90){$p$}
\put(-221,90){$p$}
\put(-360,90){$q$}
\put(-188,90){$q$}
\put(-56,135){$p$}
\put(-18,50){$q$}
\caption{This figure illustrates that $P(p,q,1)$ is in fact a 2-bridge knot. The picture in the middle
is obtained from the one the left by isotoping the arc exiting on the top left of the $p$-box. The
picture on the right is obtained by further rotating the $q$-box $180^o$ clockwise about its lower
right corner. The 2-fold branched cover of the thus obtained 2-bridge knot is the lens space $L(-pq-p-q,p+1)$.}  \label{pic6}
\end{figure}

Inspired by Lisca, we use this approach on 3-stranded pretzel knots. Thus, let
$P(p,q,r)$ denote the 3-stranded pretzel knot with $p$, $q$ and $r$ half-twists in its strands, as in
figure \ref{pic2}(b). We further assume that $p,q,r$ are odd and that
$|p|$, $|q|$, $|r| \ge 3$. In the case when any of $p$, $q$ or $r$ equals $\pm 1$, the corresponding
pretzel knot $P(p,q,r)$
is a 2-bridge knot (see figure \ref{pic6}) and so Lisca's results \cite{lisca1, lisca2} apply.

With these conventions in place, the main result of this article is the next theorem:
\begin{theorem} \label{main}
Consider the pretzel knot $P(p,q,r)$ with $p,q,r$ odd and with $|p|$, $|q|$, $|r|\ge 3$.
Then $P(p,q,r)$ is slice if and only if either
$$ p+q = 0 \quad \mbox{ or } \quad p+r=0 \quad \mbox{ or } \quad q+r=0.$$
and in each of these cases $P(p,q,r)$ is a ribbon knot. All other pretzel
knots $P(p,q,r)$ are of infinite order in the smooth knot concordance group. 
\end{theorem}
\begin{corollary}
The slice-ribbon conjecture is true for 3-stranded pretzel knots $P(p,q,r)$ with $p,q,r$ odd.
\end{corollary}
\begin{proof}
If $\min \{|p|,|q|,|r| \} \ge 3$, the corollary follows from theorem \ref{main}. If instead $\min \{|p|,|q|,|r|\} =1$,
the resulting $P(p,q,r)$ is a 2-bridge knot, in which case the corollary follows from the work of Lisca
\cite{lisca1}.
\end{proof}

\begin{remark}
When $\min\{|p|,|q|,|r|\}=1$ one direction of theorem \ref{main} still holds. Namely, if either of
$p+q$, $q+r$ or $p+r$ is zero, then $P(p,q,r)$ is still slice (see proposition \ref{ribbons} below). 
However the \lq\lq only if\rq\rq \, portion of theorem \ref{main} no longer holds: 
for example, the knot $P(23,-3,1)$ is ribbon according to \cite{lisca1}. 
\end{remark}

We mention the following easy corollary, which yields a result first proved by Ron Fintushel and
Ron Stern \cite{ronron}.
\begin{corollary} \label{ronron}
Let $P(p,q,r)$ be a pretzel knot with $p,q,r$ odd and with Alexander polynomial $\Delta_{P(p,q,r)}(t) = 1$.
Then $P(p,q,r)$ is slice if and only if $P(p,q,r)$ is the unknot.
\end{corollary}
\begin{proof}
A straightforward computation shows that the Alexander polynomial of $P(p,q,r)$ is trivial precisely when
$pq+qr+pr=-1$. If $\min \{|p|,|q|,|r|\} =  1$, say $r=1$ for concreteness,
this equation implies that $0=pq+p+q+1=(p+1)(q+1)$, showing that $p=-1$ or $q=-1$ and thus that
$P(p,q,r)$ is the unknot.

If $\min \{|p|,|q|,|r|\} \ge  3$ we can use theorem \ref{main}.
Without loss of generality, suppose that $q+r=0$.  Then the equation $pq+qr+pr=-1$ reduces to $qr=-1$, 
which is impossible for $|q|,|r| \ge 3$.
\end{proof}

Unlike for 2-bridge knots, it is not obvious that a pretzel knot $K=P(p,q,r)$ and its mirror both have
2-fold branched covers which bound negative definite plumbings. In fact, the 2-fold branched cover of $P(p,q,r)$
bounds a negative definite plumbing tree of the type considered in section \ref{theproof} (see in 
particular figure \ref{pic4}) if and only if
$\frac{1}{p}+\frac{1}{q} +\frac{1}{r} >0$ (cf. lemma \ref{negdefincidence}).  Because of this we are only able to
extract \lq\lq half\rq\rq the obstruction for sliceness that Lisca is able to use.
Using Donaldson theory we show
that there must exist some $\lambda \in
\mathbb{Z}$ such that $-q = p\lambda ^2 + r(\lambda +1)^2$
(proposition \ref{donprop}) whenever $P(p,q,r)$ is slice. It is almost magical that Heegaard Floer homology
provides the tools
which perfectly complement Donaldson's sliceness obstruction. Using the Ozsv\'ath-Szab\'o correction
terms for rational homology spheres we are further able to show that in fact only $\lambda =0,-1$
can result in a slice knot. The results about the infinite order of $P(p,q,r)$ in the 
smooth knot concordance group are proved by an extension of this argument. 

The Heegaard Floer techniques we employ suffice to prove the following
corollary which was first obtained by Andrew Casson and Cameron Gordon in
their celebrated papers \cite{cg1, cg2}.
\begin{corollary} \label{CassonGordon}
Among the twist knots $P(1,q,1)$ with $q$ odd, the only slice knots are those corresponding to
the choices $q=-1$ (the unknot) and $q=-5$ (the stevedore's knot). The only other twist knot of finite 
order in the smooth concordance group is the figure eight knot (corresponding to $q=-3$). 
\end{corollary}
Of course this corollary also follows from the work of Lisca \cite{lisca1}; we only emphasize it here
since our proof relies solely on Heegaard Floer tools.
\vskip1mm
Both sliceness obstructions which we use in the proof of theorem \ref{main} are derived from obstructions
for a 3-manifold $Y$ to bound a rational homology 4-ball $W$. In light of this, by choosing
$Y$ to the be 2-fold branched cover of the knots from theorem
\ref{main}, the latter theorem can be recast as
\begin{theorem} \label{main2}
Among the Seifert fibered spaces $M((p,1),(q,1),(r,1))$ 
with odd $p,q,r \in \mathbb{Z}$ and $|p|, |q|,|r|\ge 3$, those and only those for which either 
$p+q=0$ or $p+r=0$ or $q+r=0$ bound rational homology balls. 
\end{theorem}
In the theorem, the notation $M((\alpha _1, \beta_1),...,(\alpha _n,\beta_n))$ stands for the Seifert 
fibered space over $S^2$ obtained from $(S^2-(D_1^2\sqcup ... \sqcup D_n^2))\times S^1$ by gluing in $n$ 
solid tori $S^1 \times D^2$, the $i$-th of which is attached so that 
$$ [\partial D^2] \mapsto \alpha _i \cdot [\partial D_i^2] + \beta _i \cdot [S^1]$$
In \cite{ch} Andrew Casson and John Harer described six families of Seifert fibered spaces with three or fewer singular
fibers, all of which bound rational homology balls. As a case by case comparison shows, the examples from
theorem \ref{main2} are distinct from those appearing in the six families considered in \cite{ch}.

\vskip3mm
The rest of the paper is organized as follows: in section \ref{prelims} we review definitions
and theorems needed for the remainder of the article. Specifically, section \ref{plumbings} reviews weighted graphs and
plumbings. Section \ref{sliceness} reminds the reader of Donaldson's diagonalization theorem and explains
how it can be used to derive a sliceness obstruction for certain knots. Section \ref{sliceness2}
introduces some Heegaard Floer theory background and a second sliceness obstruction derived from correction
terms of 3-manifolds. Finally, section
\ref{algtop} which is algebro-topological in nature, is concerned with identifying which spin$^c$-structures
on a rational homology sphere extend to a rational homology ball bounded by the said sphere.
Section \ref{theproof} is devoted to the proofs of theorem \ref{main} and corollary \ref{CassonGordon}.
It is divided into two subsections, the first one explaining the input from Donaldson theory, the second
utilizing Heegaard Floer theory.

{\bf Acknowledgments } The second author would like to thank Swatee Naik and Chris Herald for stimulating 
conversations.
Additional special thanks are due to Chris Herald for his contribution of proposition \ref{ribbons} to this work.
We are indebted to Peter Ozsv\'ath for helpful email correspondence. 
\section{Preliminaries} \label{prelims}
\subsection{Weighted graphs and plumbings} \label{plumbings}
For more details on weighted graphs and the plumbing construction we refer the
reader to \cite{gompf}, Ex. 6.4.2.

\begin{figure}[htb!]
\centering
\includegraphics[width=10cm]{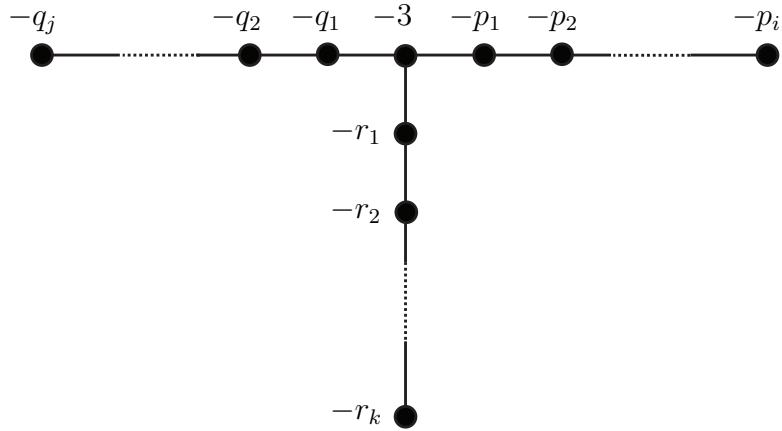}
\put(-154,153){$-3$}
\put(-125,153){$-p_1$}
\put(-96,153){$-p_2$}
\put(-18,153){$-p_i$}
\put(-185,153){$-q_1$}
\put(-215,153){$-q_2$}
\put(-292,153){$-q_j$}
\put(-170,110){$-r_1$}
\put(-170,80){$-r_2$}
\put(-170,4){$-r_k$}
\caption{The weighted graph $G=G(p,q,r)$.}  \label{pic3}
\end{figure}

Let $G$ be a finite weighted graph. Let $w(v)$ denote the weight of a vertex $v$ of $G$.
Order the vertices of $G$ in an arbitrary manner and let $v_i$ denote the $i$-th vertex.
Recall that the incidence matrix $A=A_G$ associated to $G$ with respect to such an ordering
is the matrix $A=[a_{i,j}]$ with
$$ a_{ij} = \left\{
\begin{array}{cl}
w(v_i) & \quad ; \quad \mbox{ if } i=j \cr
& \cr
\# \mbox{ of edges connecting $v_i$ and $v_j$ } & \quad ; \quad \mbox{ if } i\ne j.
\end{array}
\right.
$$
By a common abuse of notation we shall label $A_G$ simply by $G$; the context should make it clear which is meant.
Assume from now on that $G$ is a tree or a forest.
To such a weighted graph $G$ we shall associate a smooth 4-manifold with boundary $W(G)$ by plumbing
together 2-disk bundles over $S^2$ according to instructions read off from $G$. Namely, for each
vertex $v$ of $G$ we pick a disk bundle $D(v)\to S^2$ with Chern class $c_1(D(v)) = w(v)$. Given two such
disk-bundles $D(v_1)$ and $D(v_2)$, we plumb them together if and only if the vertices $v_1$ and $v_2$
are connected by an edge of $G$.  The intersection form of the resulting 4-manifold $W(G)$, expressed in terms of the basis of the spheres $S^2$ used in its construction, is the incidence matrix of $G$.
\vskip3mm
Given a pretzel knot $P(p,q,r)\subset S^3$ with $\min \{|p|,|q|,|r|\}\ge 3$, let $Y(p,q,r)$ denote the 3-manifold obtained as the
2-fold cover of $S^3$ branched along $P(p,q,r)$. These 3-manifolds are Seifert fibered spaces
with three singular fibers and their plumbing
descriptions are easily obtained according to the following recipe (see \cite{walter} or \cite{saso}).
Find continued fraction expansions of
$p/(p-1)$, $q/(q-1)$ and $r/(r-1)$:
\begin{align} \nonumber
\frac{p}{p-1} & = [p_1,p_2,...,p_i], \cr
\frac{q}{q-1} & = [q_1,q_2,...,q_j], \cr
\frac{r}{r-1} & = [r_1,r_2,...,r_k],
\end{align}
where by $[x_1,x_2,...,x_n]$ we mean
$$ [x_1,x_2,...,x_n] = x_1-\cfrac{1}{x_2-\cfrac{1}{\ddots-\cfrac{1}{x_n}}}. $$
%
Let $G=G(p,q,r)$ be the weighted graph as in figure \ref{pic3}; then $Y(p,q,r) = \partial W(G(p,q,r))$.
%
%
%
\subsection{Definite 4-manifolds and a first sliceness obstruction} \label{sliceness}
In the early 1980's, Simon Donaldson revolutionized the study of smooth 4-manifolds by introducing
Yang-Mills gauge theory as a tool for distinguishing between different smooth structures on the
same underlying topological 4-manifold. Among his most celebrated results from this time period is
the diagonalization theorem for the intersection form of a definite smooth 4-manifold,
a theorem which has become known as theorem A:
\begin{theorem}[Donaldson \cite{donaldson}]  \label{theoremA}
Let $X$ be a closed, smooth, oriented 4-manifold. If the intersection form $Q_X:H_2(X;\mathbb{Z})\otimes
H_2(X;\mathbb{Z})\to \mathbb{Z}$ is negative definite, then $Q_X$ is diagonalizable.
\end{theorem}
It is easy to see that any torsion class from $H_2(X;\mathbb{Z})$ pairs under $Q_X$ to zero with any other
class, showing that $Q_X$ descends to a map
$$Q_X: \left( H_2(X;\mathbb{Z})/\Tors \right) \otimes \left( H_2(X;\mathbb{Z})/\Tors \right)   \to \mathbb{Z}$$
Of course $H_2(X;\mathbb{Z})/\Tors$ is isomorphic to the free Abelian group $\mathbb{Z}^{b_2(X)}$.
We shall refer to the pair $(\mathbb{Z}^{b_2(X)}, Q_X)$ as a {\em lattice}.
By saying in theorem \ref{theoremA}
that \lq\lq $Q_X$ is diagonalizable\rq\rq \, we mean that the lattice $(\mathbb{Z}^{b_2(X)}, Q_X)$
is isomorphic to the standard negative definite lattice $(\mathbb{Z}^{b_2(X)}, -\mbox{Id})$ of the same
dimension.

Donaldson's theorem can be used as an obstruction to sliceness for certain knots $K\subset S^3$.
Namely, let $Y$ denote the 2-fold cover of $S^3$ branched along $K$. Suppose that $Y$ is the boundary
of a negative definite 4-manifold $W$. If $K$ is slice, then let $D^2\hookrightarrow D^4$ be a smooth
slicing disk for $K$ and let $B$ be the 2-fold cover of $D^4$ branched along the disk.
It is well known (cf. lemma 17.2 in \cite{kauffman}) 
that $B$ is a rational homology ball and $\partial B = Y$. Gluing $B$ to $W$
yields a closed, smooth 4-manifold $X$ with rk $H_2(X;\mathbb{Z}) = $ rk $H_2(W;\mathbb{Z})$.
Since, according to theorem \ref{theoremA}, the intersection form $Q_X$ is diagonalizable, 
the lattice 
$(\mathbb{Z}^{b_2(W)},Q_W)$ must embed into the negative definite lattice. The latter, as the work of Lisca \cite{lisca1} shows,
turns out to be a rather powerful obstruction for sliceness. We summarize this discussion for later use:

\begin{obstruction} \label{mainobstr}
Let $K\subset S^3$ be a knot and let $Y$ be the 2-fold cover of $S^3$ branched along $K$. Let $W$ be any
smooth
negative definite 4-manifold with $\partial W = Y$. If $K$ is slice, 
then the lattice $(\mathbb{Z}^{b_2(W)},Q_W)$ must embed in the standard negative definite intersection lattice of equal rank; that is,
there must exist a
monomorphism $\varphi :\mathbb{Z}^{b_2(W)} \to \mathbb{Z}^{b_2(W)}$ such that
$Q_W(\alpha , \beta) = -\mbox{Id}(\varphi (\alpha), \varphi (\beta))$ for any $\alpha, \beta
\in \mathbb{Z}^{b_2(W)} \cong H_2(W;\mathbb{Z})/\Tors$.
\end{obstruction}
\subsection{Heegaard Floer homology and a second sliceness obstruction} \label{sliceness2}
In a succession of seminal papers \cite{os2, os1}, Ozsv\'ath and Szab\'o constructed a
series of tools for analyzing low dimensional manifolds, a theory which has
become known as Heegaard Floer homology. In this section we focus on a particular
feature of the theory which furnishes an obstruction for a rational homology
3-sphere $Y$ to bound a rational homology ball $W$. Details
can be found in \cite{os4}, and also \cite{jabukanaik}.

The {\em correction term} $d(Y,\s)\in \mathbb{Q}$ of a spin$^c$ 3-manifold $(Y,\s)$ with
$c_1(\s)$ torsion is a rational number extracted from the Heegaard Floer homology
groups of $(Y,\s)$ \cite{os4}. 
Ozsv\'ath and Szab\'o showed that if $Y$ is a
rational homology sphere bounding a rational homology ball $W$, then $d(Y,\s)=0$
for each spin$^c$-structure $\s$ on $Y$ which extends over $W$. The correction term 
is also additive in the sense that 
\begin{equation} \label{additivityford}
d(Y_1\#Y_2, \s _1\#\s_2) = d(Y_1,\s_1) + d(Y_2,\s_2) 
\end{equation}
for any pair $(Y_1,\s_1)$ and $(Y_2,\s_2)$ of spin$^c$ rational homology $3$-spheres. 

An exercise
in algebraic topology shows that the order of $H^2(Y;\mathbb{Z})$ for $Y=\partial W$
has to be a square, say $d^2$, and that the image $\im (H^2(W;\mathbb{Z})\to H^2(Y;\mathbb{Z}))$
has order $d$. Assuming that $H^2(Y;\mathbb{Z})$ has no $2$-torsion (which is the case for the branched double-cover of a knot), we may identify $Spin^c(Y) \cong H^2(Y;\mathbb{Z})$ via the first Chern class $c_1$. Making a compatible affine identification
$Spin^c(W) \cong H^2(W;\mathbb{Z})$ (so that the diagram below commutes)
\vskip1mm
\centerline{
\xymatrix{
Spin^c(W) \ar[r]^\cong  \ar[d] & H^2(W;\mathbb{Z})  \ar[d]  \\
Spin^c(Y) \ar[r]^\cong & H^2(Y;\mathbb{Z})  \\
}
}
\vskip1mm
\noindent shows that if $Y$ bounds a rational homology ball $W$, there exists a subgroup
$V\subset H^2(Y;\mathbb{Z})$ of square root order such that $d(Y,\s)=0$
for all $\s \in V$.

When $K\subset S^3$ is a slice knot with slicing disk $D^2\hookrightarrow D^4$, we let
$W_K$ be the 2-fold cover of $D^4$ branched over the slicing disk. Then
$W_K$ is a rational homology ball and $\partial W_K=Y_K$
is just the 2-fold cover of $S^3$ branched over $K$. Using the pair $(W_K,Y_K)$ in the
above discussion yields our second sliceness obstruction:
\begin{obstruction} \label{hfhobstr}
Let $K\subset S^3$ be a knot and let $Y_K$ denote the 2-fold cover of $S^3$ branched over $K$.
If $K$ is slice, then the order of $H^2(Y_K;\mathbb{Z})$ is a square and there exists a subgroup
$V\subset H^2(Y_K;\mathbb{Z})$ of square root order such that $d(Y_K,\s)=0$ for each $\s \in V$.
The subgroup $V$ is the image of the restriction induced map $H^2(W_K;\mathbb{Z})\to H^2(Y_K;\mathbb{Z})$,
where $W_K$ is the rational homology ball obtained by a 2-fold cover of $D^4$ branched over the
slicing disk for $K$.
\end{obstruction}

To effectively utilize this obstruction, one needs to be able to
calculate the correction terms $d(Y_K,\s)$. While this is a hard problem in general,
for the case at hand we can invoke a formula due to Ozsv\'ath and Szab\'o in \cite{os6}. We summarize their results below.

Let $G$ be a negative definite weighted graph (by which we mean that its incidence
matrix is negative definite) with $n$ vertices and let $Y=Y(G) = \partial W(G)$
be the boundary of the plumbing manifold $W(G)$ associated to $G$. Identify
$H^2(W(G),Y(G);\mathbb{Z})$ with $\mathbb{Z}^n$ via the basis of Poincar\'e
duals of the $n$ attaching 2-handles, and similarly identify
$H^2(W(G);\mathbb{Z})/\Tors$ with $\mathbb{Z}^n$ via the Hom-duals of the
$n$ 2-handles.
This allows for the identification of the intersection pairing
$$H^2(W(G),Y(G);\mathbb{Z}) \otimes H^2(W(G);\mathbb{Z})/\Tors \to \mathbb{Z}$$
with the lattice $(L, G)$ with $L\cong \mathbb{Z}^n$. In addition, the long exact
sequence in cohomology for the pair $(W(G),Y(G))$, with the above bases in mind,
becomes
$$ 0 \to \mathbb{Z}^n \stackrel{G}{\to} \mathbb{Z}^n \to H^2(Y(G);\mathbb{Z})\to 0, $$
leading to the identification $H^2(Y(G);\mathbb{Z}) \cong \coker \, G$.
Let 
$G^*:L^*\otimes L^*\to \mathbb{Q}$ denote the map dual to $G$, i.e.
$G^*(G(v,\cdot), G(w,\cdot)) = G(v,w)$. 
With respect to the basis of $L^*$ dual to that of $L$, this map is expressed by the matrix $G^{-1}$.  
Recall that a covector $v \in L^*$ is called {\em characteristic} for $G$ if 
$\langle v,w \rangle \equiv G(w,w) \,(\mbox{mod} \, 2)$
for all $w\in L$. Equivalently, this is the condition that $v$, expressed in the dual basis, be congruent modulo $2$ to 
the diagonal of $G$. For a given $\s \in \coker \, G$ we shall write $Char_\s(G)$ to denote the
set of characteristic covectors whose equivalence class in $\coker \, G$ is $\s$.


A vertex $v$ of $G$ is called {\em overweight}\footnote{We prefer this more descriptive term to the use of 
{\em bad} in \cite{os6}.  Notice that if $G$ has no overweight vertices, then the corresponding knot is {\em thin} in the sense of Khovanov and knot Floer homology, further justifying our use of this term.}
if $-d(v) < w(v)$ where $w(v)$ is the weight
and $d(v)$ is the valence of $v$. With these conventions and definitions in place, we are ready
to state the formula we use for computing correction terms.
\begin{theorem}[Ozsv\'ath-Szab\'o \cite{os6}] \label{calccorrec}
Let $G$ be a negative definite forest with $n$ trees, each of 
which contains at most two overweight vertices. Let $W(G)$ be the plumbing
4-manifold associated  to $G$  and set $Y(G) = \partial W(G)$.
Then under the identification $c_1 : Spin^c(Y(G)) \to H^2(Y(G);\mathbb{Z})$, we have
%
%
\begin{equation} \label{correctermsformula}
d(Y(G),\s) = \max _{v \in Char_\s(G)}\frac{ G^{-1}(v, v) + |G|}{4}.
\end{equation}
%
\end{theorem}
This theorem is proved in \cite{os6} for the case of $n=1$. When $n>1$, it follows from
the $n=1$ case along with \eqref{additivityford}.

\subsection{Algebraic topology} \label{algtop}
This section elucidates some of the algebro-topological underpinnings invoked in the
Heegaard Floer theory portion of the proof of theorem \ref{main} in section \ref{theproof}. In particular,
the main result of this section (proposition \ref{algtopprop} below) is of crucial importance
for our application of obstruction \ref{hfhobstr}, as it explicitly identifies $V$.

To set the stage, suppose that $X^4$ is the result of attaching $n$ $2$-handles to $D^4$,
$W^4$ is a rational
homology ball, and $\del W \cong \del X \cong Y$.
The long exact sequences in cohomology of the pairs $(X,Y)$ and $(X\cup_Y W,W)$ are related by
the excision  map induced by the inclusion $(X,Y) \into (X \cup_Y W,W)$.

Observe that by our assumption on $X$ we have $H_1(X;\mathbb{Z})=0$, and thus by Poincar\'e duality
$H^3(X,Y;\mathbb{Z})=0$. The latter group is isomorphic to $H^3(X\cup_Y W,W;\mathbb{Z})$ by excision.
Furthermore, $H^1(W;\mathbb{Z})$ and $H^1(Y;\mathbb{Z})$ are trivial by the above assumptions on $W$ and $Y$.

With these observations in mind, the long exact sequences of the pairs  $(X,Y)$ and $(X\cup_Y W,W)$
simplify to
\[ \xymatrix{
0 \ar [r] & H^2(X \cup_Y W, W;\mathbb{Z})
\ar^{\alpha} [r] \ar^{i_1} [d] & H^2(X \cup_Y W;\mathbb{Z}) \ar^{\beta} [r] \ar^{i_2} [d] & H^2(W;\mathbb{Z}) \ar [r] \ar^{i_3} [d] & 0\\
0 \ar [r] & H^2(X,Y;\mathbb{Z}) \ar^{\gamma} [r] & H^2(X;\mathbb{Z}) \ar^{\delta} [r] & H^2(Y;\mathbb{Z}) \ar [r] & 0} \]
The map $i_1$ is an isomorphism by excision, whence $\im \; \gamma \subset \im \; i_2$
and the map $\delta$ embeds the quotient $(\im \; i_2) \; / \; (\im \; \gamma)$ into $H^2(Y;\mathbb{Z})$.
In addition $\im \; i_3 = \im \; i_3 \beta = \im \; \delta i_2$.  Therefore, the map $\delta$
sets up an isomorphism
$$\im \; i_3 \cong \frac{ \im \; i_2}{ \im \; \gamma}$$
Both groups
$H^2(X,Y;\mathbb{Z})$ and $H^2(X;\mathbb{Z})$ are free abelian of rank $n$ by assumption on $X$, which
implies that $H^2(X \cup_Y W, W;\mathbb{Z})$ and $H^2(X \cup_Y W;\mathbb{Z}) / \Tors$ are as well.
By abuse of notation we shall regard
$\alpha$ as a map $H^2(X \cup_Y W, W;\mathbb{Z}) \to H^2(X \cup_Y W;\mathbb{Z}) / \Tors$ and $i_2$ as a map
$H^2(X \cup_Y W;\mathbb{Z}) / \Tors \to H^2(X;\mathbb{Z})$.  Poincar\'e duality relates these maps as indicated in the
commutative diagram
\[ \xymatrix{
H^2(X \cup_Y W, W;\mathbb{Z}) \ar^{\alpha} [r] \ar^{i_1} [d] & H^2(X \cup_Y W;\mathbb{Z}) / \Tors \ar^{\varphi} [dd] \\
H^2(X,Y;\mathbb{Z}) \ar^{\psi} [d] & \\
H_2(X;\mathbb{Z}) \ar^{i_2^*} [r] & H_2(X \cup_Y W;\mathbb{Z}) / \Tors}
\]
Here $i_2^*$ denotes the map dual to $i_2$ and $\psi,\varphi$ are Poincar\'e duality maps.
It follows that we obtain the factorization
$ \gamma = i_2 \varphi^{-1} i_2^* \psi $ and hence the identification (via $\delta \circ i_2$)
\begin{equation} \nonumber
\im \; i_3 \cong \coker \; \varphi^{-1} i_2^* \psi.
\end{equation}
We summarize our findings in the next proposition.
\begin{proposition} \label{algtopprop}
Let $X$ be a 4-manifold obtained from the 4-ball by attaching 2-handles, let $W$ be a rational
homology 4-ball and assume that $\partial X \cong \partial W \cong Y$. Then the maps
\begin{align} \nonumber
i_2: & H^2(X\cup _Y W;\mathbb{Z}) \to H^2(X;\mathbb{Z}) \cr
i_2^*: & H_2(X;\mathbb{Z}) \to H_2(X\cup _Y W;\mathbb{Z})/ \Tors \cr
i_3: & H^2(W;\mathbb{Z})\to H^2(Y;\mathbb{Z})  \cr
\varphi: & H^2(X\cup _Y W;\mathbb{Z})/ \Tors \to H_2(X\cup _Y W;\mathbb{Z})/ \Tors \cr
\psi: & H^2(X,Y;\mathbb{Z}) \to H_2(X;\mathbb{Z})
\end{align}
induced by restriction, inclusion, restriction, duality, and duality, respectively, give rise to the isomorphism
\begin{equation} \label{algtop1}
\im \; i_3 \cong \coker \; \varphi^{-1} i_2^* \psi
\end{equation}
This isomorphism is facilitated by $\delta \circ i_2$ where $\delta: H^2(X;\mathbb{Z}) \to H^2(Y;\mathbb{Z})$
is the restriction induced map.
Moreover the map $\gamma :H^2(X,Y;\mathbb{Z}) \to H^2(X;\mathbb{Z})$ also induced by restriction,
can be calculated as
\begin{equation} \label{algtop2}
\gamma  = i_2 \circ \varphi ^{-1} \circ i_2^* \circ \psi
\end{equation}
\end{proposition}
\section{Proof of theorem \ref{main}} \label{theproof}
This section is devoted to proving theorem \ref{main} and corollary \ref{CassonGordon}.
The proof is divided into two
sections: the first draws on results from Donaldson theory and in particular utilizes
the sliceness obstruction \ref{mainobstr}, while the second one leans on results from 
Heegaard Floer theory and exploits obstruction \ref{hfhobstr}.

\subsection{Input from Donaldson theory} \label{DonaldsonProof}
Pretzel knots satisfy a number of symmetry relations:
\begin{equation} \label{symmetry}
P(p,q,r) = P(r,p,q)\quad \quad P(p,q,r) = P(r,q,p) \quad \quad \overline{P(p,q,r)} = P(-p,-q,-r).
\end{equation}
where $\overline{K}$ is the mirror image of $K$. Clearly $K$ is slice, ribbon, or of finite concordance order
if and only if $\overline{K}$
has the corresponding property. Thus, for convenience, and without loss of generality, we will assume that
$p$ and $r$ are positive. Remember also that the working assumption of theorem \ref{main} is that
$\min\{|p|,|q|,|r|\} \ge 3$.

When $q>0$ the signature of $P(p,q,r)$ is nonzero (lemma \ref{notnegdef} below) and so $P(p,q,r)$ cannot be 
slice. 
Thus we turn to the case of
$p,r\ge 3$ and $q\le -3$.
\begin{proposition} \label{donprop}
Consider the pretzel knot $K=P(p,q,r)$ with $p,r \ge 3$ and $q\le -3$ and all three of $p,q,r$ odd. If
$K$ is slice then there exists an integer $\lambda \in \mathbb{Z}$ such that
\begin{equation}  \label{donaldoncondition}
-q  = p \lambda ^2+ r (\lambda+1)^2
\end{equation}
\end{proposition}
%
%
%
The rest of this subsection is devoted to the proof of proposition \ref{donprop}.

The continued fraction expansions of $p/(p-1)$, $q/(q-1)$ and  $r/(r-1)$ (see section
\ref{plumbings}) are easily found to be
\begin{align} \nonumber
\frac{p}{p-1} = \overbrace{[2,2,...,2]}^{p-1} \quad  \quad \quad \frac{q}{q-1} = [1,-q+1]
\quad  \quad \quad \frac{r}{r-1} = \overbrace{[2,2,...,2]}^{r-1}
\end{align}
Therefore the plumbing graph $G=G(p,q,r)$ for $Y(p,q,r)$ (see again section \ref{plumbings})
is given as in figure \ref{pic4}(a).
\begin{figure}[htb!]
\centering
\includegraphics[width=15cm]{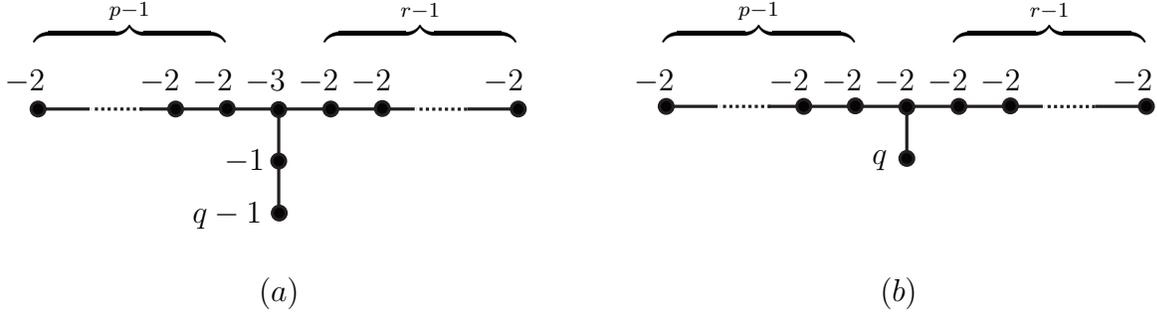}
\put(-345,50){$-3$}
\put(-325,50){$-2$}
\put(-305,50){$-2$}
\put(-255,50){$-2$}
\put(-365,50){$-2$}
\put(-385,50){$-2$}
\put(-436,50){$-2$}
\put(-353,20){$-1$}
\put(-365,0){$q-1$}
\put(-107,50){$-2$}
\put(-87,50){$-2$}
\put(-67,50){$-2$}
\put(-17,50){$-2$}
\put(-127,50){$-2$}
\put(-147,50){$-2$}
\put(-198,50){$-2$}
\put(-108,22){$q$}
\put(-340,-30){$(a)$}
\put(-105,-30){$(b)$}
\put(-425,60){$\overbrace{\phantom{iiiiiiiiiiiiiiiiii}}^{p-1}$}
\put(-315,60){$\overbrace{\phantom{iiiiiiiiiiiiiiiiii}}^{r-1}$}
\put(-187,60){$\overbrace{\phantom{iiiiiiiiiiiiiiiiii}}^{p-1}$}
\put(-77,60){$\overbrace{\phantom{iiiiiiiiiiiiiiiiii}}^{r-1}$}
\caption{The weighted graphs $G=G(p,q,r)$ (on the left) and $\widetilde G = \widetilde G (p,q,r)$
(on the right).}  \label{pic4}
\end{figure}
By blowing down the $-1$ framed vertex of $G$ we arrive at the graph $\widetilde G = \widetilde G(p,q,r)$ as in
figure \ref{pic4}(b). Since $\partial W(G) = \partial W(\widetilde G)$ we are free to use $W=W(\widetilde G)$
in obstruction \ref{mainobstr}, provided the incidence matrix of $\widetilde G$ is negative definite.
The next lemma explains when this is the case.
\begin{lemma} \label{negdefincidence}
The incidence matrix of the weighted graph $\widetilde G$ from figure \ref{pic4}(b) is negative definite
for those and only those choices of $p,q,r$ which satisfy
\begin{equation} \label{negdef}
\frac{1}{p} + \frac{1}{q} + \frac{1}{r} > 0.
\end{equation}
\end{lemma}
\begin{proof}
This follows immediately from theorem 5.2 in \cite{walter} after recognizing that the orbifold Euler characteristic
of $Y(p,q,r)$ is given by $-(\frac{1}{p} + \frac{1}{q} + \frac{1}{r})$. There is, however, a fairly easy
direct proof.  

Let $A=[a_{ij}]$ be the incidence matrix of $\widetilde G$ with respect to the following
ordering of the vertices of $\widetilde G$: Let $v_1,...,v_{p+r-1}$ be the vertices in the horizontal
branch of the graph (containing all the vertices with weight -2) with the index increasing from left to
right. Let $v_{p+r}$ be the unique vertex with weight $q$. Thus, for example, $v_p$ is the unique tri-valent
vertex.

Let $T_n=[a_{ij}]_{1\le i,j,\le n}$
be the $n\times n$ matrix in the upper left hand corner of $A$.
Recall Sylvester's criterion \cite{gilbert} from linear algebra: The matrix $A$ 
is negative definite if and only if $sign (\det T_n) = (-1)^n$ for all $n=1,...,p+r$.

For $n\le p+r-1$ the matrices $T_n$ all have the same formal shape: $-2$'s on the diagonal, 1's on the two
off diagonals and 0's everywhere else. It is easy to show (by induction on the dimension of the matrix) that
for these matrices one obtains $\det T_n  = (-1)^n (n+1)$. The determinant of $A=T_{p+r}$ is not hard to
compute either: By a last row expansion one finds that

$$\det T_{p+r} = (-1)^{p+r} ( -q(p+r) - pr) $$

Thus $A$ is negative definite if and only if $-q(p+r) -pr >0$. 
This latter condition of course translates into \eqref{negdef}.
\end{proof}
\begin{lemma} \label{notnegdef}
The signature $\sigma (P(p,q,r))$ of a pretzel knot $P(p,q,r)$ with $p,q,r$ odd, is zero if and only 
if $pq+pr+qr<0$. In particular, any pretzel knot with $pq+pr+qr \ge 0$ cannot be of finite order in the 
knot concordance group. 
\end{lemma}
\begin{proof}
The Seifert form associated to the obvious Seifert surface of $P(p,q,r)$ from 
figure \ref{pic2}(b) (i.e. the surface consisting of two disks connected by 3 bands with $p,q$ and $r$ 
half-twists respectively) shows that 
$$ \sigma (P(p,q,r)) = \sigma (S) \quad \quad \mbox{ where } \quad 
\quad S=\left[   
\begin{array}{cc}
p+q & q \cr
q & r+q
\end{array}
\right]$$ 
The matrix $S$ has signature zero if and only if its determinant $\det S = pq + pr + qr$ is negative.  
\end{proof}

With these two lemmas out of the way, we return to $W=W(\widetilde G)$. In view of lemma \ref{notnegdef},
we shall henceforth assume that $p,q,r$ satisfy condition \eqref{negdef} and thus that the incidence matrix
of $\widetilde G$ is negative definite. For simplicity we denote $\widetilde G$ by $G$ from now on.

Notice that  rk($H_2(W;\mathbb{Z}))=p+r$ and so the embedding $\varphi$ from obstruction \ref{mainobstr} is a
monomorphism
$\varphi : \mathbb{Z}^{p+r} \to \mathbb{Z}^{p+r}$. Let $\{f_1,...,f_{p+r}\}$ be the basis
of the domain of $\varphi$ which corresponds to the vertices of $G$ numbered so that
$f_{p+r}$ corresponds to the vertex with weight $q$ and $f_1,...,f_{p+r-1}$  are the vertices with weights $-2$
labeled from left to right in figure \ref{pic4}(b). With this convention the unique
trivalent vertex of $G$
corresponds to $f_p$. As a further abbreviation in notation let us write $f_i \cdot f_j$ to denote $Q_W(f_i,f_j)$.
Similarly let $\{e_1,...,e_{p+r}\}$ be a basis for the codomain of $\varphi$ and write $e_i \cdot e_j$ to denote
-Id$(e_i,e_j)$, i.e. $e_i \cdot e_j = - \delta_{ij}$.

The only linear combinations of the $e_1,...,e_{p+r}$ which yield vectors of
square $-2$ are of the form $\pm e_i \pm e_j$ for a pair of indices $i,j$.
Since $f_1$ has square $-2$, up to re-indexing of the basis $e_1,...,e_{p+r}$ and up to scaling by $-1$
we must have $\varphi (f_1) = e_1-e_2$. Since $f_2$ also has square $-2$ it too must have such a form.
However since $f_1 \cdot f_2 = 1$, one of the vectors appearing in $\varphi (f_2)$ must be either $-e_1$
or $e_2$. Thus, again up to re-indexing, we are forced to define $\varphi(f_2) = e_2-e_3$.

We proceed by induction to show that (up to a change of basis) we are forced to make the assignment
$\varphi (f_i) = e_i - e_{i+1}$ for all $i\le p+r-1$. Suppose this to be true for all $i\le m$ and consider $f_{m+1}$.
Since $\varphi (f_m) = e_m -e_{m+1}$ and since $f_m \cdot f_{m+1} = 1$, either $-e_m$ or $e_{m+1}$
is a summand of $\varphi (f_{m+1})$. If $e_{m+1}$ is a summand of $\varphi(f_{m+1})$, then since
$f_{m+1}\cdot f_i = 0$ for all $i\le m-1$, we see that the second summand of $\varphi(f_{m+1})$ cannot
be among the $e_i$, $i\le m$. Thus up to re-indexing we get $\varphi(f_{m+1}) = e_{m+1}-e_{m+2}$
completing the induction process.

On the other hand, if $-e_m$ is a summand of $\varphi(f_{m+1})$, then since $\varphi( f_{m-1} )=e_{m-1}-e_m$
and since $f_{m-1}\cdot f_{m+1} = 0$, we arrive at $\varphi (f_{m+1}) = -e_m -e_{m-1}$. But since
$\varphi ( f_{m-2} )  = e_{m-2} - e_{m-1}$, this would imply that $f_{m-2}\cdot f_{m+1} = -1$, a contradiction.
This argument fails when $m=2$ (since then $f_{m-2} = f_0$ doesn't exist). When $m=2$ there is indeed
a valid assignment of
\begin{align} \nonumber
\varphi (f_1) & = e_1 - e_2 \cr
\varphi (f_2) & = e_2 - e_3 \cr
\varphi (f_3) & = -e_1 - e_2
\end{align}
However, since $f_3 \cdot f_4=1$ we need $\varphi (f_4)$ to contain either $e_1$ or $e_2$ as a summand.
If $e_1$ were a summand of $\varphi (f_4)$ then $f_1\cdot f_4 = 0$ would imply $\varphi (f_4) = e_1 + e_2$
while if $e_2$
were a summand of $\varphi (f_4)$
then $f_2 \cdot f_4 = 0$ would imply $\varphi (f_4)  = e_2+e_3$. In either case $\varphi$
ceases to be injective, again a contradiction. Finally, observe that an $f_4$ of square $-2$ must indeed exist: Using the
symmetries of pretzel knots we have assumed that $p,r \ge 3$ and so $p+r-1\ge 5$.
Of course the number of basis elements
$f_i$ with square $-2$ is exactly $p+r-1$.

To summarize, up to a change of basis, the only possible values for
$\varphi (f_i)$ for $i=1,...,p+r-1$ are $\varphi(f_i) = e_i - e_{i+1}$.
%
%
%
%
%
%

It remains to find $\varphi (f_{p+r})= \sum _i \lambda _i e_i$. The coefficients $\lambda_i$ are readily
determined from the conditions
$$
f_j \cdot f_{p+r} = \left\{
\begin{array}{cl}
0 & ; \quad j\ne p,p+r \cr
1 & ; \quad j=p  \cr
q & ; \quad j=p+r
\end{array}
\right.
$$
The first of the equations above yields the relation $\lambda_{j+1}=\lambda_j$ for all $j\ne p,p+r$ while the
second one implies that $-\lambda_p +\lambda_{p+1} = 1$. Writing $\lambda = \lambda _p$, these two conditions
together imply that $\varphi(f_{p+r}) = \lambda (e_1+...+e_p) + (\lambda +1) (e_{p+1}+...+e_{p+r}) $.   Finally,
the condition $f_{p+r}\cdot f_{p+r} = q$ now shows that
$$ q = -p\lambda ^2 - r(\lambda+1)^2 $$
as claimed in proposition \ref{donprop}, thus completing its proof.

\subsection{Input from Heegaard Floer theory}
Recall that proposition \ref{donprop} from the previous section
imposes the relation $-q = p\lambda ^2 + r (\lambda +1)^2$ on $p,q,r$ of
any slice knot $P(p,q,r)$ (with $p,r\ge 3$ and $q\le -3$). While this is a strong restriction,
and suffices for example to yield corollary \ref{ronron}, only the values of $\lambda =-1$ (implying $p+q=0$)
and $\lambda =0$ (implying $r+q=0$) lead
to slice knots, as we shall presently see. The main inputs for the computation of this section
are the obstruction \ref{hfhobstr} and proposition \ref{algtopprop}.

For now let $K\subset S^3$ be any slice knot and let $Y_K$ be its 2-fold branched cover. Assume that
$Y_K$ bounds a negative definite plumbing $X$ associated to a weighted graph $G$ (which we assume is
a forest) with $n$ vertices. Recall that we abusively denote the incidence matrix of the graph $G$
by $G$ as well.
Let $W_K$ be the rational homology 4-ball obtained by a 2-fold cover of $D^4$ branched over the slicing
disk for $K$.

Let $\tilde f_1, ..., \tilde f_n \in H_2(X;\mathbb{Z})$ be the basis represented by the $n$ 2-handles of $X$
and let $f_1,...,f_n \in H^2(X,Y_K;\mathbb{Z})$ be the basis of their Poincar\'e duals.
Furthermore, let $e_1,...,e_n\in H^2(X;\mathbb{Z})$ be the basis of the Hom-duals of $\tilde f_i$:
$e_i(\tilde f_j) = \delta _{ij}$. With respect to these choices of bases, the restriction induced map
$\gamma:H^2(X,Y_K;\mathbb{Z})\to H^2(X;\mathbb{Z})$ is represented by the matrix $G$. The long
exact sequence of the pair $(X,Y_K)$
$$ 0 \to H^2(X,Y_K;\mathbb{Z}) \stackrel{G}{\to} H^2(X;\mathbb{Z}) \stackrel{\delta}{\to}
H^2(Y_K;\mathbb{Z}) \to 0 $$
allows us to identify $H^2(Y_K;\mathbb{Z})$ with the cokernel of $G$ (via $\delta$).
\begin{theorem} \label{generalcasetheorem}
With $K$, $Y_K$, $W_K$, $X$, $G$, $\{\tilde f_1,...,\tilde f_n\}$, $\{f_1,...,f_n\}$ and $\{e_1,...,e_n\}$ as
above, there exists a map
$H^2(X\cup _{Y_K} W_K;\mathbb{Z}) \to H^2(X;\mathbb{Z})$ whose matrix representative $A$ (with the given choices
of bases) leads to
a factorization $G=-A A^\tau $ with the additional property that (after identifying
$H^2(Y_K;\mathbb{Z})$ with $\coker \, G$) the image
$H^2(W_K;\mathbb{Z})\to H^2(Y_K;\mathbb{Z})$ is isomorphic to $ (\im \, A )/ (\im \, G)$ (via $\delta $).
\end{theorem}
\begin{proof}
The manifold $X\cup _{Y_K} W_K$ is a closed negative definite smooth 4-manifold
and so, according to Donaldson's theorem \ref{theoremA},  $(H^2(X\cup_{Y_K} W_K;\mathbb{Z}),Q_{X\cup_{Y_K} W_K})$
is isomorphic to the  standard diagonal negative definite lattice.
Let $\{a_1,...,a_n\}\subset H^2(X\cup_{Y_K} W_K;\mathbb{Z})/\Tors$ be a basis with respect to which the intersection
pairing is represented by -Id, and let $\{b_1,...,b_n\}\subset H_2(X\cup_{Y_K} W_K;\mathbb{Z})$ be its Hom-dual basis. With respect to
this pair of bases, the Poincar\'e duality map
$\varphi :H^2(X\cup_{Y_K} W_K;\mathbb{Z})/\Tors \to H_2(X\cup_{Y_K} W_K;\mathbb{Z})/\Tors$
is represented by -Id.

With the choices of bases as above, let $A$ be the $n \times n$ matrix representing the restriction induced map
$i_2:H^2(X\cup _{Y_K} W_K;\mathbb{Z}) / \Tors \to H^2(X;\mathbb{Z})$, and note that then $A^\tau$ represents the dual map
$i_2^* :H_2(X;\mathbb{Z}) \to H_2(X\cup _{Y_K}W_K;\mathbb{Z}) / \Tors$. In addition, observe that the
Poincar\'e duality map $\psi : H^2(X,Y_K;\mathbb{Z}) \to H_2(X;\mathbb{Z})$ is represented by the identity matrix.

Plugging these observations into relation \eqref{algtop2} from proposition \ref{algtopprop} immediately yields
$$ G = -A A^\tau,$$
while the isomorphism \eqref{algtop1} from the same proposition asserts that
\begin{align}\nonumber
\im (i_3:H^2(W_K;\mathbb{Z})\to H^2(Y_K;\mathbb{Z})) &  = \delta  \circ A (\coker \, (\varphi ^{-1}\circ A^\tau \circ \psi)) \cr
& = \delta \circ A \left( \frac{H^2(X\cup _{Y_K} W_K;\mathbb{Z})}{\im \, (\varphi ^{-1}\circ A^\tau \circ \psi)} \right) \cr
& = \delta \left( \frac{\im \, A}{\im (A \circ \varphi ^{-1}\circ A^\tau \circ \psi)}\right) \cr
& = \delta \left(\frac{\im \, A}{\im \, G } \right), \cr
\end{align}
where $\delta : H^2(X;\mathbb{Z}) \to H^2(Y_K;\mathbb{Z})$ is the map facilitating the identification of
$H^2(Y_K;\mathbb{Z})$ with $\coker \, G$.
\end{proof}
We now return to the case of $K=P(p,q,r)$ and $G$ the weighted graph from figure \ref{pic4}(b).
Recall that in this case we denote $Y_K$ by $Y(p,q,r)$.
The matrix $A$ whose existence is asserted by theorem \ref{generalcasetheorem} is easily determined from
the methods from section \ref{DonaldsonProof}. Namely, the embedding $\varphi :\mathbb{Z}^{p+r} \to
\mathbb{Z}^{p+r}$ from that section is precisely the adjoint of the map
$H^2(X\cup _{Y_K}  W_K;\mathbb{Z}) / \Tors \to H^2(X;\mathbb{Z})$, and thus is represented by the matrix $A^\tau$. Moreover,
the methods from section \ref{DonaldsonProof} also avouch that, up to a change of basis and a choice of
$\lambda \in \mathbb{Z}$, the matrix $A$ is uniquely determined:
$$
A = \left[
\begin{array}{rrrrrrrr}
1  & -1  & 0 & 0 & ...  & 0 & 0 & 0 \cr
0 & 1  & -1 & 0 & ...  & 0 & 0 & 0 \cr
0  & 0 & 1 & -1 & ...  & 0 & 0 & 0 \cr
0  & 0  & 0 & 1 & ...  & 0 & 0 & 0 \cr
\vdots & \vdots & \vdots & \vdots & \vdots & \vdots & \vdots & \vdots \cr
0  & 0  & 0 & 0 & ...   & -1  & 0  & 0 \cr
0  & 0  & 0 & 0 & ...  & 1  & -1  & 0\cr
0  & 0  & 0 & 0 & ...   & 0 & 1  & -1\cr
\lambda   & \lambda  & \lambda & \lambda & ...   & \lambda +1  & \lambda +1 & \lambda + 1\cr
\end{array}
\right].
$$
The $i$-th row of $A$, for $i=1,...,p+r-1$ has a $1$ in its $i$-th column, a $-1$ in its $(i+1)$-st column
and zeros elsewhere. The $(p+r)$-th row of $A$ has $\lambda$'s in its first $p$ columns and $(\lambda+1)$'s
in its remaining $r$ columns.  An easy explicit calculation shows that indeed the factorization $G=-AA^\tau$
holds.

Theorem \ref{generalcasetheorem} also tells us that the subgroup $V\subset H^2(Y(p,q,r);\mathbb{Z})$ from
obstruction \ref{hfhobstr} is isomorphic to $(\im \, A)/(\im \, G)$ via the map
$\delta : H^2(X;\mathbb{Z})\to H^2(Y(p,q,r);\mathbb{Z})$. This makes it easy to find an upper bound on
the number of
vanishing correction terms $d(Y(p,q,r),\s)$ for $\s \in V$. Towards this goal, pick $v=A x \in \im \, A$.
The term $v^\tau G^{-1} v$ from \eqref{correctermsformula} simplifies to
$$ v^\tau G^{-1} v =  - x^\tau A^\tau (A A^\tau)^{-1} A  x = - x^\tau A^\tau (A^\tau)^{-1}A^{-1} A  x = -|x|^2,$$
showing that
\begin{equation} \label{correctiontermwithx}
d(Y(p,q,r),\s) = \max_{Ax\in Char_\s(G)}  \frac{p+r-|x|^2}{4}.
\end{equation}

The requirement that $v=Ax$ be characteristic (cf. section \ref{sliceness2}) easily translates into a condition on $x$ itself:
\begin{align} \label{charx}
& v=A x \mbox{ is characteristic } \Leftrightarrow  v_i \equiv G_{ii} \,  (\mbox{mod} 2), \quad  \forall i \cr
\quad \Leftrightarrow \quad & \sum _j A_{ij}x_j \equiv  \sum _j A_{ij}^2 \quad  (\mbox{mod} 2), \quad  \forall i \cr
\quad \Leftrightarrow \quad & \sum _j A_{ij}x_j \equiv  \sum _j A_{ij} \quad  (\mbox{mod} 2), \quad  \forall i.
\end{align}
Since $\det(A)$ is odd (up to sign, it is the square-root of the knot determinant), the matrix $A$ is invertible $(\mbox{mod} 2)$, and so the vector $x (\mbox{mod} 2)$ is uniquely determined by this condition.  On the other hand, taking $x_i \equiv  1 \,  (\mbox{mod} 2), \quad   \forall i,$ clearly satisfies this equation, so it must be the unique solution. Combining \eqref{correctiontermwithx} with \eqref{charx} we see that the only way for $d(Y(p,q,r),\s)$ to be zero for a given
$\s \in V$ is that the corresponding $x=(x_1,...,x_{p+r})$ have coordinates $x_i \in \{\pm 1\}$ for all
$i=1,...,p+r$. While there are $2^{p+r}$ such vectors $v=A x$ in $\im \, A$,
there are significantly fewer equivalence classes of these vectors in
$V = (\im \, A)/(\im \, G)$.

To see this,
define $\ell : \mathbb{Z}^{p+r} \to \mathbb{Z}$ by $\ell (x) = x_1+...+x_{p+r}$. Observe that the first
$p+r-1$ columns of $A^\tau $ generate the kernel of $\ell$, showing that
any two vectors $v=A x$ and $v'=A x'$ with $x'=x+y$ and $y\in \kerr \, \ell \subset \im \, A^\tau$
belong to the same equivalence class in $V$.

Finally, the functional $\ell$, when restricted to the set $\{x\in \mathbb{Z}^{p+r} \, | \, x_i \in \{ \pm 1\} \}$,
only takes on $p+r+1$ distinct values, showing that there can be at most that many characteristic
covectors in $v\in V$ with vanishing correction terms.  But according to obstruction \ref{hfhobstr} and
the result of proposition \ref{donprop}, there need to be at least $|\cokerr(A^\tau)| = |\det(A^\tau)| = |p\lambda + r (\lambda +1)|$
vanishing correction terms. It is now a simple matter to establish that
$$ |p\lambda  + r (\lambda +1)| > p+r+1 $$
when $p,r \ge 3$, unless $\lambda = 0, -1$. We have then proved
\begin{proposition} \label{intermediate1}
Let $P(p,q,r)$ be a pretzel knot with $|p|,|q|,|r| \ge 3$ and all three $p,q,r$ odd.
If $P(p,q,r)$ is slice then either $p+q=0$ or $q+r=0$ or $p+r=0$.
\end{proposition}
Theorem \ref{main} follows from the preceding proposition along with the proceeding two. 
Proposition \ref{ribbons} was communicated to us by Chris Herald, whose
input we gratefully acknowledge (see also \cite{gompf}, p. 216).
\begin{proposition} \label{ribbons}
Any pretzel knot $P(p,q,r)$ with either $p+q=0$ or $p+r=0$ or $q+r=0$ is ribbon.
\end{proposition}
\begin{proof}
Without loss of generality (see the symmetries \eqref{symmetry}) we can assume that $q+r=0$. Recall
(cf. \cite{lisca1}) that if a
knot $K$ can be turned into an $(m+1)$-component unlink by attaching $m$ bands to it for some $m\ge 1$,
then $K$ is ribbon. A pretzel knot $P(p,q,r)$ with $q+r=0$ can easily be isotoped to a 2-component unlink
after attaching a single band as illustrated in figure \ref{pic4}.
\begin{figure}[htb!]
\centering
\includegraphics[width=13cm]{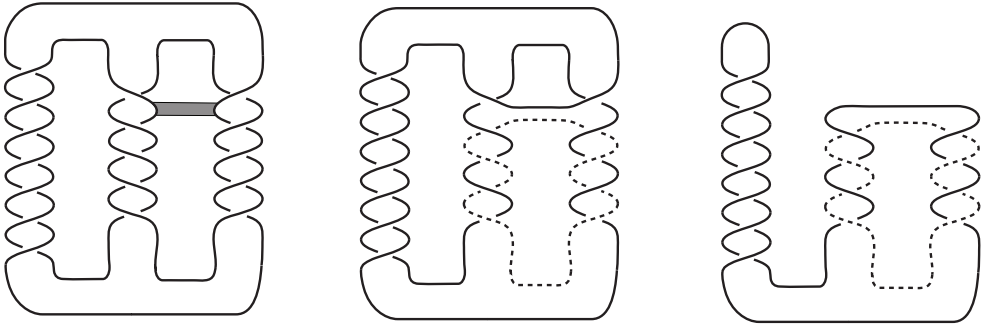}
\caption{Attaching a band to $P(7,-5,5)$ yields a 2 component unlink. }  \label{pic4}
\end{figure}
\end{proof}
\vskip3mm
\begin{proposition} \label{intermediate3}
Consider the pretzel knot $P(p,q,r)$ with $p,q,r$ odd and $|p|,|q|,|r|\ge 3$. If neither of the 
quantities $p+q$, $p+r$ and $q+r$ are zero then $P(p,q,r)$ has infinite order in the smooth 
knot concordance group. 
\end{proposition}
\begin{proof}
For convenience, let us again assume that $p,r\ge 3$ and that $q\le 3$ (when $q>0$ lemma \ref{notnegdef} shows that 
$P(p,q,r)$ cannot be of finite concordance order). 
Let $K=P(p,q,r)$ with $p,q,r$ odd and suppose that there exists some 
integer $n\ge 2$ such that $\#^n K$ (the $n$-fold connected sum of $K$ with itself) is slice. Recall that 
the branched double-cover $Y_{\#^nK}$ of $\#^nK$ is related to the branched double-cover $Y_K$ of $K$ as 
$Y_{\#^nK} = \#^n Y_K$.   Thus a plumbing description for $Y_{\#^nK}$ is given by the forest $\mathcal G$
of $n$ trees with each tree being a copy of the graph $G$ from figure \ref{pic4}(b). As before we use the
notation $\mathcal G$ to also denote the incidence matrix of this weighted graph. 

According to theorem \ref{generalcasetheorem}, there exists an $n(p+r)\times n(p+r)$ square matrix $\mathcal A$ 
such that $\mathcal{G} = - \mathcal{A} \mathcal{A}^\tau $ and the subgroup 
$V\subset \coker \, \mathcal G \cong H^2(Y_{\#^nK};\mathbb{Z})
\cong \oplus _{j=1}^n H^2(Y_K;\mathbb{Z})$ from obstruction \ref{hfhobstr} is given by $(\im \, \mathcal A)/
(\im \, \mathcal G)$. 

The matrix $\mathcal A$ is obtained as in the case of $n=1$ (see section \ref{DonaldsonProof}) although it now has $n^2$
undetermined parameters.  To see the general form of $\mathcal A$, let us pick a basis 
$\{f_1,...,f_{n(p+r)}\}$ for the lattice $(\mathbb{Z}^{n(p+r)}, \mathcal G)$ in such a way that the 
$(p+r)$-tuple $\{f_{(i-1)(p+r)+1},...,f_{i(p+r)} \}$ is to the $i$-th copy of $G$ in $\mathcal G$ as 
the basis $\{f_1,...,f_{p+r}\}$ was to the sole copy of $G$ considered in section \ref{DonaldsonProof}. In other words, 
the $f_{i(p+r)}$, $i=1,...,n$ enumerate the $n$ vertices of square $q$ in the $n$ copies of $G$ while 
for example $f_{(i-1)(p+r)+p}$, $i=1,...,n$ label all the trivalent vertices of $\mathcal G$. 
Let $\{e_1,...,e_{n(p+r)}\}$ be a basis for the standard negative definite lattice 
$(\mathbb{Z}^{n(p+r)},-\mbox{Id})$ as before and, also as before, let us label the embedding 
of $(\mathbb{Z}^{n(p+r)},\mathcal G)$ into $(\mathbb{Z}^{n(p+r)},-\mbox{Id})$ by $\varphi$. 
Such an embedding must exist according to obstruction \ref{mainobstr} and our assumption on $K$ having order $n$. 

Repeating verbatim the arguments from section \ref{DonaldsonProof} one readily finds that $\varphi(f_i) = e_i-e_{i+1}$ for 
all $i \ne (p+r), 2(p+r),...,n(p+r)$. On the other hand, writing 
$\varphi(f_{i(p+r)}) = \sum _{j=1}^{n(p+r)} \lambda _{i,j}\,  e_j$ and using the relations 
\begin{equation} \label{multirelations}
 f_{i(p+r)} \cdot f_{k} = 
\left\{
\begin{array}{cl}
q                    & \quad ; \quad k = i(p+r) \cr
1                    & \quad ; \quad k = (i-1)(p+r)+p \cr
0                    & \quad ; \quad \mbox{otherwise } \cr 
\end{array}
\right.
\end{equation}
quickly leads to the equations 
$$
\begin{array}{rcl}
k\ne i  & \quad \Longrightarrow \quad & \lambda _{i,(k-1)(p+r)+1}  = \lambda _{i,(k-1)(p+r)+2} = ... =
\lambda _{i,k(p+r)}   = : b_{i,k}\cr
& & \cr
k=i & \quad \Longrightarrow \quad & \lambda _{i,(i-1)(p+r)+1} = \lambda _{i,(i-1)(p+r)+2} = ... = 
\lambda _{i,(i-1)(p+r)+p}:= b_{i,i}  \cr
& & \lambda _{i,(i-1)(p+r)+p+1}= \lambda _{i,(i-1)(p+r)+p+2}=...=\lambda _{i,i(p+r)} \cr
& & \lambda _{i,(i-1)(p+r)+p+1} = \lambda _{i,(i-1)(p+r)+p} + 1 .
\end{array} 
$$
Substituting these back into the expression for $\varphi (f_{i(p+r)})$ yields 
\begin{equation} \label{formulaforfipplusr}
\varphi (f_{i(p+r)}) = \left( \sum _{j=1}^n b_{i,j} \left( e_{(j-1)(p+r)+1} + ... + e_{j(p+r}\right) \right) + 
\left( e_{(i-1)(p+r)+p+1} + ...+ e_{i(p+r)}\right)
\end{equation}
for any $i=1,...,n$.  
Let us assemble the $n^2$ integers $b_{i,j}$ into a matrix $B=[b_{i,k}]$. We then see that $\mathcal A$ 
is, up to choosing $B\in M_{n}(\mathbb{Z})$, the unique matrix whose 
$i$-th row is made up by the coefficients appearing in the expansion of $\varphi (f_i)$ 
in terms of the basis $\{e_1,..,e_{n(p+r)}\}$. Choosing $k= (p+r),2(p+r),...,n(p+r)$ in 
equation \eqref{multirelations}, shows that $B$ is constrained by the matrix equation 
\begin{equation} \label{matrixeq} 
- q \cdot I = p \cdot B B^\tau + r (B+I)(B+I)^\tau, 
\end{equation}
(compare this to the result of proposition \ref{donprop} in the case of $n=1$) but is otherwise arbitrary.  

Proceeding in analogy with the case of $n=1$, we define $n$ linear functionals 
$\ell _i :\mathbb{Z}^{n(p+r)} \to \mathbb{Z}$, $i=1,...n$ as 
$\ell _i(x_1,...,x_{n(p+r)}) = \sum _{j=1}^n x_{(i-1)(p+r)+j}$ and assemble these to obtain 
$\ell =(\ell_1,...,\ell_n):\mathbb{Z}^{n(p+r)} \to \mathbb{Z}^n$.


By repeating verbatim the 
argument from the case of $n=1$, it is easy to see that $v= \mathcal A x$ is characteristic for $\mathcal G$ 
if and only if all $x_i$, $i=1,...,n(p+r)$ are odd, and moreover the correction term 
$d(Y_{\#^nK},[v])$ vanishes precisely when $x_i \in \{\pm 1\}$ for each $i=1,...,n(p+r)$ (here 
$[\cdot ]$ indicated the equivalence class of $v\in \coker \, \mathcal G \cong H^2(Y_{\#^nK};\mathbb{Z})$). 
As before, the count of equivalence classes $[\mathcal A x] \in \coker \, \mathcal G$ 
equals the count of equivalence classes $[x]\in \coker \, \mathcal A^\tau$. 

Observe that the kernel of $\ell _i$ is spanned by the columns of $\mathcal A^\tau$ ranging 
from the $[(i-1)(p+r)+1]$--st to the $[i(p+r)-1]$--st. Consequently, any two $x,x'\in \mathbb{Z}^{n(p+r)}$ 
with $x-x'\in \ker \ell$ induce the same equivalence class in $\coker \, \mathcal A^\tau$. As $\ell$ 
restricted to the set $\mathcal S = \{x\in \mathbb{Z}^{n(p+r)} \, | \, x_i \in \{\pm 1\}, i=1,...,n(p+r) \}$ 
only attains $(p+r+1)^n$  different values, there can be at most that many vanishing correction terms
for $Y_{\#^n K}$. But obstruction \ref{hfhobstr} dictates that there be at least 
$\sqrt{|H^2(Y_{\#^nK};\mathbb{Z})|} = |pq+qr+pr|^{n/2}$ vanishing correction terms, showing that $K$ cannot be 
of finite order unless 
\begin{equation} \label{fijame}
\sqrt{|pq+qr+pr|} \le (p+r+1).
\end{equation}
Assuming \eqref{fijame} we proceed. First off, note that \eqref{fijame} implies that 
$$ - q \le \frac{(p+r+1)^2+pr}{p+r}= p+r +2 +\frac{1+pr}{p+r}< p+r+2+\min \{ p ,r\}, $$ 
where the strict inequality uses the assumption $p,r\ge 3$ from the beginning of the proof.  This inequality 
for $-q$ forces $b_{i,i}\in \{0,-1\}$ for all $i=1,...,n$, for otherwise the $(i,i)$--th entry from equation 
\eqref{matrixeq} yields 
\begin{align} \nonumber
 -q b_{i,i} & = p \left( b_{i,1}^2+...+b_{i,n}^2\right ) + 
 r \left( b_{i,1}^2+ ... + b_{i,i-1}^2+(b_{i,i}+1)^2+ b_{i,i+1}^2+...+b_{i,n}^2\right) \cr
 & \ge p\,  b_{i,i}^2 + r (b_{i,i} + 1)^2 \cr
 & \ge p+r + 3 \min \{p,r\} \cr
 & > p+r+2+\min\{p,r\},
\end{align} 
which contradicts equation \eqref{fijame}. Similarly, if either  
$|b_{i,j}|>1$ for some $i\ne j$ or 
if $|b_{i,j_1}|=1=|b_{i,j_2}|$ for at least two indices $j_1, j_2\ne i$, then equation \eqref{matrixeq}
again implies that $-q>p+r+2+\min\{p,q\}$, another contradiction.  Lastly, if no $j \ne i$ exists so that $b_{i,j} \ne 0$, then $-q$ equals $p$ or $r$, which is the case of a ribbon knot that we have already handled.  
%
%

We are therefore left to consider the case where, for each $i=1,...,n$, there exists exactly one $j=j(i)\ne i$ 
for which $b_{i,j(i)} = \pm 1$ and where $b_{i,i} \in \{0,-1\}$. In this case, $-q = p+r+ \min \{ p,r \}$.  The mod 2 version of equation \eqref{matrixeq}
is $0 \equiv B+B^\tau (\mbox{mod } 2)$ and accordingly $b_{i,j}\ne 0$ if and only if 
$b_{j,i}\ne 0$. This observation further limits the form of $B$. For example, consider 
the pair $\varphi (f_{p+r})$ and $\varphi(f_{j(p+r)})$ with $j=j(1)$:
\begin{align} \nonumber
\varphi (f_{p+r}) & = b_{1,1}(e_1+...+e_p)+(b_{1,1}+1)(e_{p+1}+...+e_{p+r}) + 
b_{1,j}(e_{(j-1)(p+r)+1}+...+e_{j(p+r)}) \cr
\varphi (f_{j(p+r)}) & = b_{j,1}(e_1+...+e_{p+r})+ \cr 
& \quad + b_{j,j}(e_{(j-1)(p+r)+1}+...+e_{(j-1)(p+r)+p})+(b_{j,j}+1)(e_{(j-1)(p+r)+p+1}+...+e_{j(p+r)}). 
\end{align}
By re-ordering of our bases, we can without loss of generality assume that $j=2$. 
Since $f_{p+r}\cdot f_{2(p+r)} = 0$, the above shows that 
$$b_{2,1}(p\, b_{1,1}+r(b_{1,1}+1)) + b_{1,2}(p\, b_{2,2} +  r\, (b_{2,2} +1)) = 0. $$
Choosing various values  for $b_{1,1}, b_{2,2}\in \{ 0,-1\}$ leads to the following four possibilities
(up to a further change of bases):
\begin{equation} \label{blocks}
\left[ 
\begin{array}{rr}
b_{1,1} & b_{1,2} \cr
b_{2,1} & b_{2,2}
\end{array}
\right] = 
\left[ 
\begin{array}{rr}
0 & 1 \cr
-1 & 0
\end{array}
\right], \, 
\left[ 
\begin{array}{rr}
-1 & 1 \cr
-1 & -1
\end{array}
\right], \, 
\left[ 
\begin{array}{rr}
0 & 1 \cr
1 & -1
\end{array}
\right] \, \mbox{ or }  \, 
\left[ 
\begin{array}{rr}
0 & -1 \cr
-1 & -1
\end{array}
\right].
\end{equation}
If $r<p$ only the first case can occur, if $p<r$ only the second, while if $p=r$, all four cases 
may happen. Repeating this argument for other pairs of indices 
(and we can always assume, after change of basis, that $j(i) = i\pm 1$) we see that $B$ is a block direct sum 
$B=B_1\oplus B_2\oplus ... \oplus B_{n/2}$ of $(n/2)$ $2\times 2$ matrices, with each $B_i$ 
one of the the cases on the  right-hand side of \eqref{blocks}. As a curiosity, note that
this already excludes all odd $n\ge 3$ from being the concordance order of $P(p,q,r)$. 

With this in place, we are ready to obtain a sharper bound on the number of classes in $\coker(\mathcal A^\tau)$ 
represented by $\mathcal S$, and use this bound to conclude the argument in our remaining case.  
The decomposition of $B$ implies that $\mathcal A^\tau$ decomposes as a direct sum of $2(p+r) \times 2(p+r)$ 
matrices, and $\coker(\mathcal A^\tau)$ decomposes accordingly as the direct sum of cokernels 
of these matrices.   

For concreteness, let us assume that $B_1$ equals the first $2\times 2$ matrix on the right-hand side of 
\eqref{blocks}, and let $\mathcal A_1$ 
denote the corresponding $2(p+r)\times 2(p+r)$ summand of $\mathcal A$.  In this case, we have 
\begin{align} \nonumber
\varphi (f_{p+r}) & = \left( e_{p+1}+\cdots+e_{p+r}\right) + \left( e_{p+r+1}+\cdots+e_{2(p+r)}\right) \cr
\varphi( f_{2(p+r)}) & = -(e_1+ \cdots +e_{p+r})+(e_{p+r+p+1}+\cdots+e_{2(p+r)}) \cr
\varphi (f_j) & = e_j-e_{j+1} \quad \quad \mbox{ for all } j\ne p+r, 2(p+r), j\le 2(p+r),
\end{align}
and so $\mathcal A_1^\tau $ takes on the form 
$$
\mathcal A_1^\tau  = 
\left[
\begin{array}{rrrrr|rrrrr}
1   & 0   & ... & 0 & 0 & 0 & 0 & ...& 0 & -1 \cr
-1  & 1   & ... & 0 & 0 & 0 & 0 & ...& 0 & -1 \cr
0  & -1   & ... & 0 & 0 & 0 & 0 & ...& 0 & -1 \cr
\vdots & \vdots & \vdots & \vdots & \vdots & \vdots & \vdots & ...& \vdots & \vdots \cr
0 & 0& ... & 1 & 1 & 0 & 0 & ...& 0 & -1 \cr
0 & 0& ... & -1 & 1 & 0 & 0 & ...& 0 & -1 \cr \hline 
0 & 0& ... & 0 & 1 & 1 & 0 & ...& 0 & 0 \cr
0 & 0& ... & 0 & 1 & -1 & 1 & ...& 0 & 0 \cr
0 & 0& ... & 0 & 1 & 0 & -1 & ...& 0 & 0 \cr
\vdots & \vdots & \vdots & \vdots & \vdots & \vdots & \vdots & ...& \vdots & \vdots \cr
0 & 0& ... & 0 & 1 & 0 & 0 & ...& 1 & 1 \cr
0 & 0& ... & 0 & 1 & 0 & 0 & ...& -1 & 1 \cr
\end{array}
\right].
$$
Thus the $j$-th column of $\mathcal A_1^\tau$, with the exception of $j=p+r, 2(p+r)$, contains a $1$ 
in the $j$-th row and a $-1$ in the $(j+1)$-st row and zeros in all other rows. 
The $(p+r)$-th column has zeros in its first $p$ rows, 
followed by $p+2r$ $1$'s. The $2(p+r)$-th column has $-1$'s in its first $p+r$ rows, followed by $p$ zeros 
which in turn are followed by $r$ $1$'s.

The column vectors of $\mathcal A_1^\tau$, excluding the $(p+r)$-th and the $2(p+r)$-th column,  
form a basis for the kernel of the map 
$\ell_{1,2} := (\ell_1,\ell_2) : \mathbb{Z}^{2(p+r)} \to \mathbb Z^2$.
Therefore, any two vectors $x,y \in \mathbb{Z}^{2(p+r)}$ with $x-y \in \ker \, \ell _{1,2}$ induce the 
same class in $\coker \, \mathcal A_1^\tau $. Since $\ell _{1,2}$, when restricted to the set 
$\mathcal S_1:= \{ x\in \mathbb{Z}^{2(p+r)}\, |\, x_i \in \{\pm 1\} \}$ takes on only $(p+r+1)^2$ 
values, there can be at most this many vanishing correction terms associated to this first block $B_1$ of $B$. 
These $(p+r+1)^2$ values can be understood concretely: The equality $\ell _{1,2} (x) = (k_1,k_2)$ means that 
$x$ has $k_1$ copies of $-1$ in its first $p+r$ coordinates and $k_2$ copies of $-1$ in its second $p+r$ 
coordinates, here $k_i \in \{0,...,p+r\}$. Thus we can think of vectors $[x] \in \coker \, \mathcal A_1^\tau $ 
with vanishing correction terms as being enumerated by the set $\{ (k_1,k_2) \, | \, k_i\in \{0,...,p+r\}\}$. 

Taking the $(p+r)$-th and the $2(p+r)$-th columns of $\mathcal A_1^\tau$ into account, we see that there are 
further 
equivalences among such pairs $(k_1,k_2)$, namely  
$$
\begin{array}{rlc} 
(k_1,p+r) & \sim (k_1-r, 0) \quad \quad & \mbox{ for } \quad k_1 = r,r+1,...,r+p, \cr
(0,k_2) & \sim (p+r,k_2-r)  \quad \quad & \mbox{ for } \quad k_2 = r,r+1,...,r+p.
\end{array}
$$
These equivalences then further reduce the number of $[x]\in \coker \, \mathcal A_1^\tau$ whose associated correction
terms vanish, to at most $(p+r+1)^2 - 2(p+1)$. Since $\mathcal A$ contains $n/2$ such blocks, there are 
at most $[(p+r+1)^2-2(p+1)]^{n/2}$ vanishing correction terms, yet according to 
obstruction \ref{hfhobstr} we need to have at least $(-q(p+r) -pr)^{n/2} = ((p+2r)(p+r)-pr)^{n/2}$ 
many.
Finally, the inequality %
$$  (p+2r)(p+r)-pr \le (p+r+1)^2-2(p+1)$$
simplifies to $(r-1)^2\le 0$, leading to the contradiction $r=1$ (since we assumed $\min\{p,r\}\ge 3$). 
The remaining three possibilities for $B_i$ are treated in the same vein and lead to either $r=1$ or $p=1$, 
both contradictions, thus completing the proof of proposition \ref{intermediate3}.
\end{proof}

Observe that taking $r=1$ in the case $-q=p+r+\min \{ p,r \}$ treated above leads to the family of knots $P(p,-(p+2),1)$.  Each knot in this family is amphicheiral, hence has order two in the smooth concordance group.  This fact provides some justification for the intricacy involved in handling this case.

\vskip5mm

Theorem \ref{main} follows from propositions \ref{intermediate1} -- \ref{intermediate3}.  
\vskip5mm
\noindent {\em Proof of corollary \ref{CassonGordon} } Notice that a necessary condition for sliceness of
the twist knot $K=P(1,q,1)$ is that its determinant be a  square:
$$\det K =  -2q-1 = \ell ^2$$
for some odd integer $\ell$. The knots $P(1,q,1)$ with $q>0$ have nonzero signature (according to lemma 
\ref{notnegdef}) and can
consequently not be slice. Therefore assume $q<0$ henceforth.
The double branched cover $Y(1,q,1)$ of
$P(1,q,1)$ is the boundary of the plumbing 4-manifold $X$ associated to the
weighted graph $G$ with only a pair of vertices with weights $-2$ and $q$ and connected by a single edge.
The incidence matrix
$G$ can be factored as $G=-AA^\tau$ as shown:
$$ G= \left[
\begin{array}{cc}
-2 & 1 \cr
1 & q
\end{array}
\right]
\quad \quad \mbox{ and } \quad \quad
A= \left[
\begin{array}{cc}
1 & -1 \cr
\lambda & \lambda +1
\end{array}
\right],
$$
for some $\lambda$ with $\lambda ^2+(\lambda +1)^2=-q$.
We would like to note that while this factorization
follows from theorem \ref{generalcasetheorem} and the condition on $\lambda$ from
proposition \ref{donprop} (whose proofs rely
on Donaldson's theorem \ref{theoremA}), these facts can easily obtained directly by just solving the matrix
equation
$G=-AA^\tau$ for $A$.

Since $H^2(Y(1,q,1);\mathbb{Z})$ is cyclic of order $\ell ^2$, it has a unique
subgroup $V$ of order $\ell$. On the other hand, the order  of $(\im \, A)/(\im \, G)$ is $\ell$, forcing $V=(\im \, A)/(\im \, G)$. We would like to remark once more that while this statement too
follows from theorem \ref{generalcasetheorem}, in the present case we can obtain it \lq\lq by hand\rq\rq
without having to rely on Donaldson's theorem.

A covector $v= A x$ is characteristic for $G$ precisely
when $x_i\equiv 1 \, (\mbox{mod} \, 2)$. Arguing as in the proof of theorem \ref{main}, we see from this
that there
are at most 3 vanishing correction terms of $Y(1,q,1)$, while obstruction \ref{hfhobstr} dictates that there
need to be at least $\ell$ of them. Thus we get $\ell \le 3$ and hence $\ell =1$ or $\ell = 3$. Since
\begin{align} \nonumber
\ell & = 1 \quad \Longrightarrow \quad q=-1 \quad \Longrightarrow \quad  P(1,q,1) = \mbox{ unknot }, \cr
\ell & = 3 \quad \Longrightarrow \quad q=-5 \quad \Longrightarrow \quad P(1,q,1) = \mbox{ the stevedore's knot },
\end{align}
these values indeed lead to slice knots. 

Finally, equation \eqref{fijame} from the proof of theorem \ref{main} was derived without the use 
of the assumption $\min\{ |p|,|q|,|r|\} \ge 3$ and therefore remains valid in the present context. It implies 
that any finite order twist knot $P(1,q,1)$ must obey the inequality 
$$ -2q-1\le 9.$$
The only other $q$ satisfying this relation (besides $q=-1,-5$) is $q=-3$, as claimed.

\hfill $\square$
\vskip3mm
As a concluding remark, we would like to point out that our proof of the Fintushel-Stern result
(corollary \ref{ronron}) only relied on Donaldson's diagonalization theorem. It is not hard, however, to
come up with a purely Heegaard Floer theoretic proof of their result: it follows from work of
C. Livingston \cite{chuck} and
Ozsv\'ath-Szab\'o \cite{os13} that the concordance invariant $\tau$ (cf. \cite{os11}) for three stranded
pretzel knots $P(p,q,r)$ with $p,q,r$ odd and $p,r >0$ is given by
$$ \tau (P(p,q,r)) = \left\{
\begin{array}{cc}
-1 & \quad ; \quad -q <\min \{ p,r \} \cr
& \cr
0 & \quad ; \quad -q \ge \min \{ p,r \}
\end{array}.
\right.
$$
Therefore, any slice pretzel knots $K=P(p,q,r)$ with trivial Alexander polynomial
$\Delta_K(t) = 1$  satisfies the inequality $-q \ge \min\{p,r\}$ and the equation
$pq+pr+qr = -1$. For concreteness assume that $r\le p$. Then $-q\ge r$ implies that
$$ 1 =  -q(p+r) -pr \ge r(p+r)-pr = r^2 \quad \Longrightarrow \quad r=1.$$
But then
$$ 0 = 1 + pr + qr + pq = 1+p+q+pq = (1+p)(1+q), $$
from which we infer that $q=-1$ (since $p>0$). Clearly $P(p,-1,1)$ is the unknot.

\bibliographystyle{plain}

\end{document}